\documentclass[12pt,leqno]{article}
\usepackage{amsmath,amssymb,amscd,latexsym,bm}

\usepackage[all]{xy}

\newcommand{\C}{{\mathbb{C}}}

\newcommand{\Ge}{\mathbb{G}}

\newcommand{\Q}{{\mathbb{Q}}}
\newcommand{\oQ}{\overline{\Q}}
\newcommand{\R}{{\mathbb{R}}}

\newcommand{\Z}{{\mathbb{Z}}}

\newcommand{\hZ}{\hat{\Z}}

\newcommand{\tf}{\tilde{f}}

\newcommand{\Bott}{\mathrm{Bott}}

\newcommand{\ch}{\mathrm{ch}\,}

\newcommand{\ddet}{\mathrm{det}}

\newcommand{\id}{\mathrm{id}}

\renewcommand{\mod}{\mathrm{mod}\,}

\newcommand{\rel}{\mathrm{rel}}

\newcommand{\spec}{\mathrm{spec}\,}

\newcommand{\GL}{\mathrm{GL}}

\newcommand{\Imm}{\mathrm{Im}\,}

\newcommand{\Ker}{\mathrm{Ker}\,}

\newcommand{\RRe}{\mathrm{Re}}

\newcommand{\Sing}{\mathrm{Sing}}

\newcommand{\Tr}{\mathrm{Tr}}
\newcommand{\ttop}{\mathrm{top}}

\newcommand{\Ah}{{\mathcal A}}
\newcommand{\Dh}{\mathcal{D}}

\newcommand{\Hh}{\mathcal{H}}

\newcommand{\Nh}{{\mathcal N}}
\newcommand{\Oh}{{\mathcal O}}

\newcommand{\onA}{\overline{A}}

\newcommand{\tDelta}{\tilde{\Delta}}

\newcommand{\dsigma}{\dot{\sigma}}

\newcommand{\dalpha}{\dot{\alpha}}

\newcommand{\hotimes}{\hat{\otimes}}

\newcommand{\iso}{\xrightarrow{\; \sim \;}}

\newcommand{\colim}{\mathrm{colim}}

\newtheorem{theorem}{Theorem}[section]

\newtheorem{theodef}[theorem]{Theorem / Definition}

\newtheorem{cor}[theorem]{Corollary}

\newcounter{aufz}
\newenvironment{rem}{\noindent {\bf Remark}}{}

\newenvironment{definition}{\noindent {\bf Definition}\it}{}

\newenvironment{proof}{\noindent {\bf Proof}}{\mbox{}\hspace*{\fill}$\Box$}

\newcommand{\Fh}{\mathcal{F}}

\newcommand{\verk}{\mbox{\scriptsize $\,\circ\,$}}

\newcommand{\Nbb}{\mathbb{N}}

\newcommand{\Rbb}{\mathbb{R}}

\newcommand{\Zbb}{\mathbb{Z}}

\newcommand{\sgn}{\operatorname{sgn}}

\newcommand{\sHom}{\mathcal{H}\kern -.5pt om}
\newcommand{\sExt}{\mathcal{E}\kern -.5pt xt}
\begin{document}
\title{Entropy, Mahler measure, operator determinants and Chern characters}
\author{Christopher Deninger\footnote{University of M\"unster, Mathematical Institute, Einsteinstr. 62, 48149 M\"unster, Germany, c.deninger@uni-muenster.de; Funded by the Deutsche Forschungsgemeinschaft (DFG, German Research Foundation) under Germany's Excellence Strategy EXC 2044--390685587, Mathematics M\"unster: Dynamics--Geometry--Structure and the CRC 1442 Geometry: Deformations and Rigidity} \and Thu H\`a Trieu\footnote{Hanoi University of Science and Technology, Faculty of Mathematics and Informatics, 1 Dai Co Viet, Bach Mai Ward, Hanoi, Vietnam, ha.trieuthu@hust.edu.vn; Funded by DAAD Research Grant 57746808}}
\date{}
\maketitle
\begin{abstract}
We give a much deeper conceptual explanation than previously known for certain a priori surprising relations between higher Deligne-Beilinson regulators on the one hand and entropies, Mahler measures and Fuglede-Kadison determinants on the other hand.

\end{abstract}
\section{Introduction} \label{sec:1}
For a Laurent polynomial $f$ in $R = \Z [t^{\pm 1}_1 , \ldots , t^{\pm 1}_d]$ the number
\[
m (f) = \int_{T^d} \log f \, d \mu 
\]
is called the (logarithmic) Mahler measure of $f$. Here $T^d = (S^1)^d$ is the real $d$-torus and $\mu$ is its Haar probability measure. For an algebraic number $\alpha \in \oQ$ with minimal polynomial $f$, the quantity $m (f)$ is the unnormalized Weil height of $\alpha$. See \cite{B} for a nice survey why Mahler measures are relevant for diophantine approximation and for Lehmer's problem. In \cite{LSW} it is shown that $m (f)$ is the entropy of a natural $\Z^d$-action on a compact abelian group naturally attached to $f$. Since the 1980's starting with \cite{S}, several polynomials $f$ have been found for which $m (f)$ can be expressed in terms of zeta- and $L$-values from number theory. For example
\begin{equation}
\label{eq:1}
m (3 + t_1 + t_2) = \frac{2 \sqrt{3}}{3 \pi} L (\chi , 2) + \log 3 \; ,
\end{equation}
where $\chi$ is the non-trivial character of $(\Z / 3)^{\times}$, so that
\[
L (\chi , s) = 1 - \frac{1}{2^s} + \frac{1}{4^s} - \frac{1}{5^s} + \frac{1}{7^s} - \frac{1}{8^s} \pm \cdots
\]
For Laurent polynomials $f$ which do not vanish in any point of $T^d$, such relations were explained to some extent in \cite{D1} by the Beilinson conjectures on special values of $L$-functions and the formula
\begin{equation}
\label{eq:2}
m (f) = \langle r_{\Dh} \{ f , t_1 , \ldots , t_d \} , (2 \pi i)^{-d} [T^d] \rangle \; .
\end{equation}
Here $\{ f , t_1 , \ldots , t_d \}$ is an element (a symbol) in the algebraic $K$-group $K_{d+1} (R_f)$ of $R_f = R [1 / f]$ and $r_{\Dh} \{ f , t_1 , \ldots , t_d \}$ is its image in the real Deligne cohomology under the Deligne-Beilinson regulator map $r_{\Dh}$. In our situation the Deligne cohomology group is canonically isomorphic to singular cohomology in degree $d$ with coefficients $\R (d) = (2 \pi i)^d \R$ of the space 
\[
U = \{ z \in (\C^{\times})^d \mid f (z) \neq 0 \} \; .
\]
In \eqref{eq:2} we evaluate on the homology class $(2 \pi i)^{-d} [T^d]$ of $T^d \subset U$ (recall that by assumption $f$ does not vanish on $T^d$). Now one is not far away from the situation where the Beilinson conjectures imply that in very favourable cases an evaluation as in \eqref{eq:2} is essentially an $L$-value at an integer. For example, this principle suggested the first relation between the Mahler measure of a two-variable polynomial and the special value $L(E,2)$ for an elliptic curve $E/\mathbb{Q}$; see \cite{D1}. In fact, numerical evidence for such relations was provided for a family of two-variable polynomials in \cite{Boy98}. To date, these relations have been rigorously established for only finitely many polynomials in this family (see \cite{Bru16}, \cite{Lal10}, \cite{LSZ16}, \cite{LR07}, \cite{RV99}, \cite{RZ12}, \cite{RZ14}, \cite{Zud14}). In 2006, further numerical conjectures relating the Mahler measures of three-variable polynomials to special values $L(E,3)$ of the $L$-functions of elliptic curves were proposed in \cite{Boy06}; see also \cite{Bru20}. These conjectural relations were subsequently investigated in the framework of Beilinson's conjectures; see, for example, \cite{Lal15}, \cite{Tri26}, \cite{Bru20K4}. Using the method of modular regulators developed in \cite{BZ23}, the first instance of such a relation between the Mahler measure of a three-variable polynomial and $L(E,3)$ was proved in \cite{Bru23}. This confirmed a consequence of the Beilinson conjectures whose general proof remains elusive. A relation between the Mahler measure of a polynomial in four variables and the special value $L(X, 4)$ for a singular $K3$ surface $X$ was first studied in \cite{Tri25}.

While the proof of \eqref{eq:2} given in \cite[Proposition 1.2]{D1} is a straightforeward calculation, the first author always wondered about the deeper meaning of formula \eqref{eq:2}. At the time of writing the paper \cite{D1}, the only conceptual interpretation of $m (f)$ for $d \ge 2$ was as an entropy. But for what reason should entropy be related to a regulator map on higher $K$-theory?  

In \cite[Example 3.13]{L} an interpretation of $m (f)$ as the logarithm of the Fuglede-Kadison determinant of $f$ was given, where $f$ is viewed as a unit in the von Neumann algebra $\Nh = L^{\infty} (T^d , \mu)$ of the group $\Gamma = \Z^d$, see section \ref{sec:2}. In fact, the relation of entropy of ``algebraic'' actions with Fuglede-Kadison determinants turned out to be very close and to persist for all countable amenable and even sofic groups $\Gamma$ and $f \in \Z \Gamma$, c.f. \cite{D2}, \cite{DS}, \cite{LT}, \cite{Hay}. 

The question about the conceptual meaning of formula \eqref{eq:2} now becomes: Why should the Fuglede-Kadison determinant of $f$ be related to a Deligne-Beilinson regulator on higher algebraic $K$-theory? The second author suggested that the work of Tamme \cite{T1}, \cite{Tam12} relating the Deligne Beilinson regulator with Karoubi's relative Chern character should be relevant. Eventually this insight led to a complete understanding of the situation and of formula \eqref{eq:2} as follows.

The original definition of the Fuglede-Kadison determinant is functional-analytic but in \cite{HS} a beautiful interpretation via algebraic and topological $K$-theory in low degrees was found. See also \cite{Har}. This new description has a natural reformulation in terms of Karoubi's relative $K$-theory of topological algebras developed in \cite{Kar87} and in an improved form in \cite{Kar97}. We review this material in section \ref{sec:3}.

Using these theories the first natural contact of the Mahler measure Fuglede-Kadison determinant with $K$-theory is established, more precisely with Karoubi's relative Chern character in degree one from $K^{\rel}_1 (\Nh)$ to the $0$-th cyclic homology group $HC_0 (\Nh)$. 

In order to pass from $K^{\rel}_1$ to higher $K$-theory, for technical reasons we have to work with the Fr\'echet algebra $C^{\infty} (T^d)$ instead of the von Neumann algebra $L^{\infty} (T^d , \mu)$. As a consequence we have to assume that $f \, |_{T^d} \in K_1 (C^{\infty} (T^d))$ is the image of an element $\tf \in K^{\rel}_1 (C^{\infty} (T^d))$. Multiplying $f$ by a uniquely determined monomial, this can always be achieved. 

The assembly map for the group $\Gamma = \Z^d$ produces a natural class $\{ T^d \}$ in $K_d (R)$ and hence in $K_d (C^{\infty} (T^d))$ from the fundamental class in $H_d (T^d , \Z) = H_d (\Gamma , \Z)$. Explicitly we have $\{ T^d \} = \{ t_1 , \ldots , t_d \}$, the ordering $t_1 , \ldots , t_d$ corresponding to the orientation. This is the intrinsic, coordinate free meaning of the symbol $\{ t_1 , \ldots , t_d \}$ in formula \eqref{eq:2}. 

Using the module structure of $K^{\rel}_*$ under $K_*$, the product of $\tf$ with the fundamental class $\{ T^d \} = \{ t_1 , \ldots , t_d \}$ gives a class in $K^{\rel}_{d+1} (C^{\infty} (T^d))$. Its image under Karoubi's relative Chern character $ch^{\rel}_{d+1}$ gives back the Mahler measure of $f$ as one easily sees using equivariance of $ch^{\rel}_{d+1}$. We formalize these considerations by introducing higher analogues of the de la Harpe Skandalis determinant. 

Finally we compare Karoubi's relative Chern character for the Fr\'echet algebra $C^{\infty} (T^d)$ with the Deligne-Beilinson regulator map for the affine variety $\spec (R_f \otimes \C)$ using several results of Tamme's thesis \cite{T1}, \cite{Tam12}. 

The various steps that we just described explain formula \eqref{eq:2} conceptually if we view $m (f)$ as the logarithmic Fuglede-Kadison determinant of $f$.

We end this introduction by mentioning three further research directions. 

In the present paper, we dealt exclusively with the case where the Laurent polynomial $f$ does not vanish in any point of $T^d$ or equivalently where $f$ is a unit in the von Neumann algebra $\Nh$. Extending our approach to the general case poses many interesting questions.

In \cite{D3}, \cite{D4} theories of $p$-adic entropy and $p$-adic Fuglede-Kadison determinants were developed and a relation with a regulator from algebraic $K$-theory to syntomic cohomology was found. Karoubi has generalized his theory to $p$-adic topological algebras and Tamme has made the comparison of Karoubi's $p$-adic relative Chern character with the syntomic regulator at least for proper varieties \cite{T2}. So, as in the real case, there should be a natural explanation in the $p$-adic case for the formulas in \cite[3.3]{D4}.

In a different direction, it is a challenging problem to generalize the passage from $K^{\rel}_1$ to higher $K^{\rel}_i$ for \textit{non-commutative} groups $\Gamma$ and to interpret the resulting formula. The problem is that the algebraic $K$-theory of non-commutative algebras has no ring structure. However, we have ideas in this direction and will explore them in the future. This plan is one reason why we described the background theory in more generality than the case $\Gamma = \Z^d$ would require. The other reason is to lend credibility to our view on the Mahler measure by embedding it into a much bigger conceptual background. 

\textit{Acknowledgement} The authors are grateful to Shanshan Hua, Jens Kaad, Max Karoubi, Wolfgang L\"uck, Thomas Nikolaus and Maxime Ramzi for discussions and email exchanges at the formative stage of this work. The first author would also like to thank Pedro F. dos Santos very much for the invitation to the Istituto superior t\'ecnico where some work on the project was done. Ai use was restricted to free models as occasional search engines. 
\section{Mahler measure, entropy and Fuglede-Kadison determinants} \label{sec:2}

In this section we give a review of the topics in the title. This will clarify the conceptual meaning of the Mahler measure integral and embed it into a conceptual context. 

Consider a countable discrete group $\Gamma$. Its von~Neumann algebra $\Nh \Gamma$ has the following description. Let $L^2 (\Gamma)$ be the Hilbert space of square summable complex sequences $(x_{\gamma})_{\gamma \in \Gamma}$ which we view as formal series $\sum_{\gamma} x_{\gamma} \gamma$. The group $\Gamma$ acts by left multiplication on such series and this gives a unitary representation of $\Gamma$ on $L^2 (\Gamma)$. The von~Neumann algebra of $\Gamma$ is the sub $*$-algebra of $\Gamma$-equivariant bounded operators on $L^2 (\Gamma)$, thus $\Nh \Gamma = B (L^2 (\Gamma))^{\Gamma}$. The unitary representation $\Gamma \to B (L^2 (\Gamma)) , \gamma \mapsto R_{\gamma^{-1}}$, where $R_{\gamma}$ denotes right multiplication by $\gamma$, commutes with the left $\Gamma$-action and hence defines a homomorphism $\Gamma \to (\Nh \Gamma)^{\times}$. It therefore induces a $\C$-algebra homomorphism $\C \Gamma \to \Nh \Gamma$ which by evaluating on the unit $e \in \Gamma$ is easily seen to be injective. We will view $\C \Gamma$ as a subalgebra of $\Nh \Gamma$ via this inclusion. It can be shown that $\Nh \Gamma$ is the strong closure of $\C \Gamma$ in $B (L^2 (\Gamma))$. The norm-continuous functional $\tau_{\Gamma} : \Nh \Gamma \to \C$ defined by $\tau_{\Gamma} (A) = (A (e), e)$ is a trace in that $\tau_{\Gamma} (AB) = \tau_{\Gamma} (BA)$. It is positive in the sense that $\tau_\Gamma (AA^*) \ge 0$ for all $A \in \Nh \Gamma$. Restricted to $\C \Gamma$ we have $\tau_{\Gamma} (\sum_{\gamma} a_{\gamma} \gamma) = a_e$. For a unit $A \in (\Nh \Gamma)^{\times}$ the Fuglede-Kadison determinant is defined by the functional calculus: The operator $AA^*$ is positive and $\log AA^* \in B (L^2 (\Gamma))$ is $\Gamma$-equivariant, so that $\log AA^* \in \Nh \Gamma$. Set
\begin{equation}
\label{eq:3}
\ddet_{\Nh \Gamma} (A) = \exp \Big( \frac{1}{2} \tau_{\Gamma} (\log AA^*) \Big) \; .
\end{equation}
The main result of \cite{FK} asserts that
\[
\ddet_{\Nh \Gamma} : (\Nh \Gamma)^{\times} \longrightarrow \R^{> 0}
\]
is a \textit{homomorphism} with $\ddet_{\Nh\Gamma} (A^*) = \ddet_{\Nh\Gamma} (A)$. 

Morally, the reason for this is the Baker-Campbell-Hausdorff formula, but the actual proof is different. The inclusion $\C \Gamma \subset \Nh \Gamma$ extends by continuity to an inclusion $L^1 (\Gamma) \subset \Nh \Gamma$ of the $L^1$-group algebra of $\Gamma$. In particular, units in $L^1 (\Gamma)$ give units in $\Nh \Gamma$. Following \cite[Example 3.13]{L} we now recall how Mahler measures of Laurent polynomials fit into this context. For the group $\Gamma = \Z^d$, we may identify $L^2 (\Gamma)$ with $L^2 (T^d)$ via the Fourier transform $\Fh$. Under the resulting identification of $B (L^2 (\Gamma))$ with $B (L^2 (T^d))$, the $\Gamma$-equivariant operators on $L^2 (\Gamma)$ become the multiplication operators by essentially bounded functions on $T^d$. This give a $*$-isomorphism
\begin{equation}
\label{eq:4}
\Nh \Gamma \iso L^{\infty} (T^d , \mu) \; , \; A \mapsto \Fh (A (0)) \;,
\end{equation}
here $0 \in \Z^d \subset L^2(\Z^d)$.  Complex Laurent polynomials in $d$ variables may be viewed as functions on $T^d$. The resulting $\C$-algebra homomorphism
\[
\C [t^{\pm 1}_1 , \ldots , t^{\pm 1}_d] \subset L^{\infty} (T^d , \mu) = \Nh (\Z^d)
\]
is injective since $\C \Gamma \subset \Nh \Gamma$. Let us recall \cite[Proposition 3.1]{D2}: For an element $f \in \C [\Z^d]$ the following conditions are equivalent:\\
a) $f$ is a unit in $L^1 (\Z^d)$,\\
b) $f$ is a unit in $\Nh (\Z^d)$,\\
c) $f$ viewed as a Laurent polynomial does not vanish in any point of $T^d$.

Under these conditions we have
\begin{equation}
\label{eq:5}
\log \ddet_{\Nh\Gamma} (f) = \int_{T^d} \log |f| \, d\mu = m (f) \; .
\end{equation}
The only non-trivial implication c) $\Rightarrow$ a) follows from a theorem of Wiener asserting that $L^1 (\Gamma)^{\times} = L^1 (\Gamma) \cap (\Nh \Gamma)^{\times}$ for $\Gamma = \Z^d$. Formula \eqref{eq:5} follows immediately from \eqref{eq:3} and \eqref{eq:4}.

The other context where Mahler measures have a conceptual meaning is as entropies in the theory of algebraic dynamical systems as mentioned in the introduction. For a countable discrete group $\Gamma$ and $f \in \Z \Gamma$ the discrete abelian group $\Z \Gamma / (\Z \Gamma) f$ has a $\Gamma$-action by group automorphisms via left multiplication. This gives a $\Gamma$-action on the compact abelian Pontrjagin dual $X_f = \widehat{\Z \Gamma / (\Z \Gamma) f}$ by continuous group automorphisms which necessarily preserve the Haar measure. By definition, $X_f$ is a closed subshift of $\widehat{\Z \Gamma} = (\bigoplus_{\Gamma} \Z)^{\wedge} = \prod_{\Gamma} \hZ = (S^1)^{\Gamma}$. An action of $\Gamma$ on a compact metrizable space $X$ is called expansive if for some (equivalently all) metrics $d$ on $X$ there exists some $\varepsilon > 0$ such that for all $x \neq y$ in $X$ there is some $\gamma \in \Gamma$ with $d (\gamma x , \gamma y) \ge \varepsilon$. By \cite[Theorem 3.2]{DS} $X_f$ is expansive if and only if $f \in L^1 (\Gamma)^{\times}$. For the $\Gamma$-action on $\Gamma = \Z^d$ this is equivalent to $f (z) \neq 0$ for all $z \in T^d$, by Wiener's theorem. Since the  $\Gamma$-action on $X_f$ preserves Haar-measure, the entropy $h_f$ of the action is defined if $\Gamma$ is amenable \cite{D2} or even if it is only sofic \cite{Hay}. For expansive actions of residually finite groups $\Gamma$ it is relatively straightforeward to prove the formula
\begin{equation}
\label{eq:6}
h_f = \log \ddet_{\Nh\Gamma} f \quad \text{and in particular} \; h_f = \int_{T^d} \log |f| \, d \mu \; \text{for} \; \Gamma = \Z^d \; .
\end{equation}
Namely, in the expansive case, $h_f$ can be expressed as a renormalized limit of logarithms of numbers of periodic points. Correspondingly $\ddet_{\Nh\Gamma} f$ is the limit of renormalized determinants on $\C \Gamma^{(n)}$ where the groups $\Gamma^{(n)}$ are finite quotients of $\Gamma$. The respective finite terms agree and hence also the limits. See \cite{DS} for details. Formula \eqref{eq:6} is also true for non-expansive actions under mild conditions on $f$ but we will not be concerned with this case which even for $\Gamma = \Z^d$ is much more complicated, \cite{LSW}, \cite{LT}, \cite{Hay}.
\section{Karoubi's relative Chern character} \label{sec:3}
In this section we review a few results from \cite{Kar97} which improve on \cite{Kar87}. Using the completed projective tensor product $\hotimes_{\pi}$ instead of the algebraic tensor product one defines for Fr\'echet algebras $A$, Hochschild homology $HH_* (A)$, cyclic homology $HC_* (A)$, negative cyclic homology $HC^-_* (A)$ and periodic cyclic homology $HP_* (A)$. These are possibly non-Hausdorff topological vector spaces which fit into a commutative diagram of continuous linear maps with exact rows
\begin{equation}
\label{eq:7}
\vcenter{
\xymatrix{
\ldots \ar[r] & HC_{i-1} (A) \ar[r]^{\partial} \ar@{=}[d] & HC^-_i (A) \ar[r]^I \ar[d] & HP_i (A) \ar[r]^S \ar[d]_{\lambda} & HC_{i-2} (A) \ar[r]^-{\partial} \ar@{=}[d] & \ldots \\
\ldots \ar[r] & HC_{i-1} (A) \ar[r]^{\partial} & HH_i (A) \ar[r]^I & HC_i (A) \ar[r]^S & HC_{i-2} (A) \ar[r]^-{\partial} & \ldots
}}
\end{equation}
Here $S$ is Connes' periodicity operator. We now recall the definitions of these homologies in analogy with the treatment of the algebraic case in \cite[section 2]{BO}. For two $\Z$-graded vector spaces $V_{\bullet}$ and $W_{\bullet}$ define the completed graded tensor product by the formula
\begin{equation}
\label{eq:8}
(V_{\bullet}  \hotimes W_{\bullet})_n = \prod_{p+q=n} V_p \otimes W_q \; .
\end{equation}
Below the $V_p$ will be finite dimensional. For an associative, unital Fr\'echet algebra $A$ over $k = \R$ or $\C$ set $\onA = A / k$ and consider the (normalized, topological) Hochschild complex $C_{\bullet} (A)$ where $C_n (A) = A \hotimes_{\pi} (\onA)^{\hotimes_{\pi} n}$ for $n \ge 0$ and $C_n (A) = 0$ for $n < 0$, see \cite[section 2.11]{Kar87}. The Hochschild- and Connes boundary maps $b : C_n (A) \to C_{n-1} (A)$ resp. $B : C_n (A) \to C_{n+1} (A)$ are first defined on elementary tensors by the following formulas and then extended continuously using the universal property of the completed projective tensor product:
\begin{align*}
b (a_0 \otimes \ldots \otimes a_n) & = \sum^{n-1}_{i=0} (-1)^i a_0 \otimes \ldots \otimes a_i a_{i+1} \otimes \ldots \otimes a_n + (-1)^n a_n a_0 \otimes a_1 \otimes \ldots \otimes a_{n-1}\\
B (a_0 \otimes \ldots \otimes a_n) & = \sum^n_{i=0} (-1)^{ni} 1 \otimes a_i \otimes \ldots \otimes a_n \otimes a_0 \otimes \ldots \otimes a_{i-1} \footnotemark \;,
\end{align*}
\footnotetext{Notice that the map $B$ mentioned in \cite[section 2.3]{Kar87} is defined for the unormalized Hochschild complex, and the second sum in the formula actually vanishes after passing to the normalized Hochschild complex.}see \cite[sections 2.2, 2.3]{Kar87}. The formulas $b^2 = 0 , B^2 = 0$ and $bB + Bb = 0$ hold. We view $k [u]$ and $k [u , u^{-1}]$ as graded $k$-algebras with $\deg u = -2$ and $\deg u^{-1} = 2$. Using the tensor product \eqref{eq:8} we define complexes of Fr\'echet spaces by setting\footnote{These complexes are defined analogously to those in Karoubi's work. In particular, the cyclic complex $CC_\bullet$ is the normalized topological counterpart of the unnormalized cyclic complex $\mathrm{Tot}(B_{**})$ introduced in \cite[Section 2.3]{Kar87}. The formal variable $u$, of degree $-2$, is used to keep track of the column index in the bicomplex.}
\begin{align*}
CC_{\bullet} (A) & = k [u , u^{-1}] / uk [u] \otimes C_{\bullet} (A), \\
CC^-_{\bullet} (A) & = k [u] \hotimes C_{\bullet} (A), \\
CP_{\bullet} (A) & = k [u, u^{-1}] \hotimes C_{\bullet} (A),
\end{align*}
with differential given by $\id \otimes b + u (\id \otimes B)$. Note that for $CC_{\bullet} (A)$, the tensor products $\otimes$ and $\hotimes$ give the same complex. Note also, that countable products of Fr\'echet spaces are Fr\'echet spaces. The homologies $HH_* (A) , HC_* (A) , HC^-_* (A)$ and $HP_* (A)$ are defined as the homologies of the complexes $C_{\bullet} (A) , CC_{\bullet} (A) , CC^-_{\bullet} (A)$ and $CP_{\bullet} (A)$. Consider the commutative diagram of graded $k$-vector spaces with exact rows
\begin{equation}
\label{eq:9}
\vcenter{
\xymatrix{
0 \ar[r] & k [u] \ar[r] \ar[d] & k [u, u^{-1}] \ar[r]^-{\cdot u} \ar[d] & k [u , u^{-1}] / uk [u] \ar[r] \ar@{=}[d] & 0 \\
0 \ar[r] & k \ar[r] & k [u, u^{-1}] / uk [u] \ar[r]^-{\cdot u} & k [u , u^{-1}] / uk [u] \ar[r] & 0 \, .
}}
\end{equation}
Tensoring this diagram with $C_{\bullet} (A)$ using $\hotimes$ of \eqref{eq:8} we obtain a commutative diagram of complexes of Fr\'echet spaces
\begin{equation}
\label{eq:10}
\vcenter{
\xymatrix{
0 \ar[r] & CC^-_{\bullet} (A) \ar[r] \ar[d] & CP_{\bullet} (A) \ar[r]^{u \otimes \id} \ar[d] & CC_{\bullet} (A) [2] \ar[r] \ar@{=}[d] & 0 \\
0 \ar[r] & C_{\bullet} (A) \ar[r] & CC_{\bullet} (A) \ar[r]^{u \otimes \id} & CC_{\bullet} (A) [2] \ar[r] & 0 \; .
}}
\end{equation}
Here $K_{\bullet} [a]$ is the usual $a$-shifted complex $K_{\bullet}$. The long exact homology sequence of \eqref{eq:10} gives the diagram \eqref{eq:7}. The Connes operator $S$ comes from $u$-multiplication, $I$ from the inclusion, $\partial$ is the boundary map, and $\lambda$ is induced from the projection.

For the Fr\'echet algebra $A = C^{\infty} (M)$ of smooth $k = \R$ or $\C$-valued functions on a closed smooth manifold $M$ the various homology groups can be calculated as follows. Let $\Ah^{\bullet} (M)$ denote the complex of smooth $k$-valued forms on $M$. The unique continuous maps
\[
\pi_n : C_n (A) \longrightarrow \Ah^n (M)
\]
with
\[
a_0 \otimes \ldots \otimes a_n \longmapsto \frac{1}{n!} a_0 da_1 \wedge \ldots \wedge d a_n
\]
intertwine $b$ with $0$ and $B$ with $d$. By the topological version of the Hochschild-Kostant-Rosenberg theorem due to Connes, $\pi$ induces a quasi-isomorphism
\[
\pi : (C_{\bullet} (A) , b) \longrightarrow (\Ah^{\bullet} (M) , 0) \; .
\]
Thus $HH_{\bullet} (A) = \Ah^{\bullet} (M)$. The following induced maps of complexes are also quasi-isomorphisms, where the differential on the right is $D = u (\id \otimes d)$:
\begin{align*}
CC_{\bullet} (A) & \longrightarrow k [u , u^{-1}] / uk [u] \otimes \Ah^{\bullet} (M), \\
CC^-_{\bullet} (A) & \longrightarrow k [u] \hotimes \Ah^{\bullet} (M), \\
CP_{\bullet} (A) & \longrightarrow k [u, u^{-1}] \hotimes \Ah^{\bullet} (M) \; .
\end{align*}
Let us write down explicitly the cyclic homology for the algebra of smooth functions on $M$. An element $\alpha \in k [u , u^{-1}] / uk [u] \otimes \Ah^n (M)$ can be written as a finite sum 
\[\alpha = \sum_{[-n/2] \le i \le 0} u^{i} \otimes \omega_i, \quad \omega_i\in \Ah^{n+2i}(M).\]
We have $D(\alpha) = \sum_{[-n/2] \le i \le -1} u^{i+1} \otimes d\omega_i$  as the term for $i=0$ vanishes in the quotient $k[u, u^{-1}]/uk[u]$. Hence, $\alpha$ is a cycle if and only if $\omega_i$ is closed for every $i \le -1$ and $\omega_0 \in \Ah^n  (M)$ is arbitrary. The cycles of degree $n$ are therefore given by
\[
Z_n = \Ah^n(M) \oplus u^{-1} Z^{n-2}(\Ah^\bullet, d) \oplus u^{-2} Z^{n-4}(\Ah^\bullet, d) \oplus \cdots.
\] 
The boundaries of degree $n$ of $k [u , u^{-1}] / uk [u] \otimes \Ah^{\bullet}(M)$ are given by
\[
B_n  = d \Ah^{n-1} (M) \oplus u^{-1} d \Ah^{n-3}(M) \oplus u^{-2} d \Ah^{n-5}(M) \oplus \cdots
\]  
This give an ismorphism 
\[
HC_n (A) = \Ah^n (M) / d \Ah^{n-1} (M) \oplus \bigoplus_{i \le -1} H^{n+2i}_{dR} (M , k).
\]
Similarly, there are isomorphisms 
\begin{align*}
HC^-_n (A) & = \Ker (d : \Ah^n (M) \to \Ah^{n+1} (M)) \oplus \bigoplus_{i \ge 1} H^{n+2i}_{dR} (M,k), \\
HP_n (A) & = \bigoplus_{i\in \Z} H^{n+2i}_{dR} (M,k) \; .
\end{align*}
If one substitutes these descriptions of the cyclic homologies into diagram \eqref{eq:7}, the map $\partial$ is induced by $d$ and the other maps are the obvious ones.

The $0$-th topological $K$-group of $A$ is by definition equal to its $0$-th algebraic $K$-group defined in terms of finitely generated projective $A$-modules, $K^{\ttop}_0 (A) = K_0 (A)$. Let $\GL (A) = \GL_{\infty} (A)$ be the colimit of the groups $\GL_r (A)$ via the inclusions $\GL_r (A) \hookrightarrow \GL_{r+1} (A) , M \mapsto \left( \begin{smallmatrix} M & 0 \\ 0 & 1 \end{smallmatrix} \right)$. 
The higher algebraic $K$-groups of $A$ are obtained from the simplicial set $B_{\bullet} \GL (A) = (\GL (A)^n)_{n \ge 0}$ with the face and degeneracy operators by multiplication respectively insertion of the unit at the appropriate places, 
\[
K_i (A) = \pi_i (|B_{\bullet} \GL (A)|^+) \quad \text{for} \; i \ge 1 \; .
\]
Here $|\;|$ is the geometric realization and $+$ denotes the Quillen plus-construction with respect to the (perfect) commutator subgroup of $\pi_1 (|B_{\bullet} \GL (A)|) = \GL (A)$. By definition $K_i (A)$ depends only on $A^{\delta}$, the ring $A$ with the discrete topology. For Banach algebras, or more generally for Fr\'echet algebras $A$ whose units are open in $A$, the higher topological $K$-groups are defined by the formula
\begin{equation}\label{KtopBanach}
	K^{\ttop}_i (A) = \pi_{i-1} (\GL (A) , 1) \quad \text{for} \; i \ge 1 \; .
\end{equation}
For general Fr\'echet algebras, Karoubi (see e.g. \cite[section 3.1]{CK}) uses a different definition of $K^{\ttop}_i (A)$. Consider the standard $n$-simplex
\begin{equation}\label{stardardsimplex}
	\Delta^n := \{(x_0, \dots, x_n) \in \Rbb^{n+1} \ | \ x_i \ge 0, \sum_{i=0}^n x_i = 1\}.
\end{equation}  
Let $A_\bullet = (A_n)_{n \ge 0}$ be the simplicial algebra with $A_n = C^\infty(\Delta^n, A)$ the algebra of smooth functions on $\Delta^n$ with values in $A$. By a smooth function on $\Delta^n$, we mean the  restriction of a smooth function with compact support in $\Rbb^{n+1}$.  For $r \in \Nbb \cup \{\infty\}$, consider the simplicial set $\GL_r(A_\bullet) = (\GL_r(A_n))_{n\ge 0}$. Let $B_\bullet \GL_r(A_\bullet)$ be the classifying simplicial set of $\GL_r(A_\bullet)$, that is, the diagonal of the bisimplicial set $[p], [q] \mapsto B_p \GL_r(A_q)$. 

\begin{definition} (Karoubi) \label{Ktop} The higher \textit{topological $K$-groups} of a Fr\'echet algebra $A$ are 
	\[K_i^\mathrm{top}(A)  := \pi_i(B_\bullet \GL(A_\bullet)), \quad i \ge 1 \; .\] 
\end{definition}

\begin{rem}
If the units of $A$ are open in $A$, it is known that $\GL_r (A_n) = C^{\infty} (\Delta^n , \GL_r (A))$, so that with the smooth and continuous singular chain complexes we have: 
\[
|\GL_r (A_{\bullet})| = |\Sing_{C^{\infty}} (\Delta^{\bullet} , \GL_r (A))| \simeq |\Sing (\Delta^{\bullet} , \GL_r (A)| \simeq \GL_r (A) \; .
\]
Hence, in this case, Karoubi's definition of $K^{\ttop}_i (A)$ agrees with the previous one mentioned in \eqref{KtopBanach}. For general Fr\'echet algebras, however, the inverse of a smooth map $\Delta^n \to \GL (A)$ is not necessarily smooth so that $\GL_r (A_n)$ may be strictly smaller than $C^{\infty} (\Delta^n , \GL_r (A))$. For dense Fr\'echet algebras $A$ in a $C^*$-algebra $B$, such that $A$ is stable under the holomorphic functional calculus in $B$, the units of $A$ are open in $A$. Moreover in this case $K^{\ttop}_i (A) = K^{\ttop}_i (B)$ for all $i \ge 0$. This applies in particular to the Fr\'echet algebra of smooth functions on a compact manifold $M$, where $A = C^{\infty} (M) \subset B = C^0 (M)$. It follows that
\[
K^{\ttop}_i (C^{\infty} (M)) = K^{\ttop}_i (C^0 (M)) = K^{-i}_{\ttop} (M) \; .
\]
\end{rem}

In \cite{Tam12} Tamme suggested the following improved definition of Karoubi's relative $K$-groups \cite[3.1]{CK}, \cite[\S\,5]{Kar87} which has good functorial properties.



Consider the integral completion functor $\Zbb_\infty$ defined by Bousfield and Kan on the category of simplicial sets. It has the following main properties (see \cite[Ch. 1, \S\,5]{BK72}).
	\begin{enumerate}
		\item [(i)] A map $f : S_\bullet \to T_\bullet$ induces an isomorphism on integral homology if and only if $\Zbb_\infty (f) : \Zbb_\infty S_\bullet \to \Zbb_\infty T_\bullet$ is a weak homotopy equivalence.
		\item [(ii)] If $S_\bullet$ has the homotopy type of an $H$-space, the natural morphism $S_\bullet \to \Zbb_\infty S_\bullet$ is a weak homotopy equivalence.
	\end{enumerate}

Now let us consider the following diagram 
\begin{equation*}
	\xymatrix{B_\bullet \GL(A) \ar[rr]^f \ar[d] && |B_\bullet \GL(A)|^+ \ar[d]\\
	\Zbb_\infty B_\bullet \GL(A) \ar[rr]^{\Zbb_\infty (f)} && \Zbb_{\infty} |B_\bullet \GL(A)|^+}
\end{equation*}
By Quillen's plus-construction, $f$ induces an isomorphism on integral homology, hence $\Zbb_\infty (f)$ is a weak homotopy equivalence. As $|B_\bullet \GL(A)|^+$ is an $H$-space, the map $|B_\bullet \GL(A)|^+ \to \Zbb_\infty |B_\bullet \GL(A)|^+$ is also a weak homotopy equivalence. This implies Dror's theorem 
\[
K_i(A) = \pi_i (|B_\bullet \GL(A)|^+) \simeq \pi_i(\Zbb_\infty B_\bullet \GL(A)).
\]
Since $B_\bullet \GL(A_\bullet)$ has the homotopy type of an $H$-space, the natural map $B_\bullet \GL(A_\bullet) \to \Zbb_\infty B_\bullet \GL(A_\bullet)$ is a weak homotopy equivalence. By definition, we have 
\[K_i^\mathrm{top}(A) = \pi_i(\Zbb_\infty B_\bullet \GL(A_\bullet)).\]
The natural map from the algebraic $K$-group $K_i(A)$ to the topological $K$-group $K_i^\mathrm{top}(A)$ is induced by the map $\Zbb_\infty B_\bullet \GL(A) \to \Zbb_\infty B_\bullet \GL(A_\bullet).$ 

\begin{definition}(Tamme)
	Karoubi's higher \textit{relative $K$-groups} of a Fr\'echet algebra  $A$ are defined by 
	\[
	K_i^\mathrm{rel}(A) := \pi_i (F),\quad i \ge 1,
	\]
	where $F$ is the homotopy fiber of the map $\Zbb_\infty B_\bullet \GL(A) \to \Zbb_\infty B_\bullet \GL(A_\bullet).$
\end{definition} 

From the definition we obtain the following exact sequence
\begin{equation}\label{karoubi seq}
	\cdots \to K_{i+1}(A) \to K_{i+1}^\mathrm{top}(A) \to K_i^\mathrm{rel}(A) \to K_i(A) \to K_i^\mathrm{top}(A)\to \cdots.
\end{equation}

Karoubi defines a relative Chern character $\ch^{\rel}_i$ which fits into the following commutative diagram with exact rows where $\ch^-_i$ and $\ch^{\ttop}_i$ are Chern characters on algebraic and topological $K$-theory, \cite{Kar97}, \cite{Kar87}
\begin{equation}\label{diag1}
\vcenter{
\xymatrix{\cdots \ar[r] & K_i^\mathrm{rel}(A) \ar[r] \ar[d]^{\mathrm{ch}^\mathrm{rel}_i} & K_i(A) \ar[r] \ar[d]^{\mathrm{ch_i^{-}}} & K_i^\mathrm{top}(A)\ar[d]^{\mathrm{ch}_i^\mathrm{per}} \ar[r] & K_{i-1}^\mathrm{rel}(A) \ar[r]\ar[d]^{\mathrm{ch}_{i-1}^\mathrm{rel}} &\cdots \\
		\cdots \ar[r] & HC_{i-1}(A)\ar[r]^{\partial} & HC^-_i(A)\ar[r]^I & HP_{i}(A) \ar[r]^S & HC_{i-2}(A) \ar[r]&\cdots.}}
\end{equation}

For Banach algebras $A$, we now recall the explicit construction of the Karoubi relative Chern character used in \cite{Kaa11}. For each integer $n \ge 0$, we normalize the standard $n$-simplex as follows:
\begin{equation*}
	\Delta^n = \begin{cases}
		\{(t_1, \dots, t_n) \in \Rbb^n\ | \ t_i \ge 0, \sum_{i=1}^n t_i \le 1\} & n \ge 1,\\
		\{0\} &  n = 0.
	\end{cases}
\end{equation*}
The vertices $v_0, \dots, v_n \in \Delta^n$ are given by 
\begin{equation*}
	v_i = \begin{cases}
		(0, \dots, 1, \dots, 0) & i \in \{1, \dots, n\},\\
		(0, \dots, 0) & i = 0.
	\end{cases}
\end{equation*} 
The reason for this presentation is  to emphasize the role of the zeroth vertex $v_0 \in \Delta^n$.  For each $r \in \Nbb \cup \{\infty\}$, we consider the simplicial set $R_r(A)_\bullet = (R_r(A)_n)_{n\ge 0}$, where $R_r(A)_n$ is the set of normalized smooth maps 
\begin{equation*}
	\sigma : \Delta^n \to \GL_r(A), \quad \sigma(v_0) = 1_r.
\end{equation*}
The face and degeneracy operators are given by 
\begin{equation*}
	\begin{aligned}
	\partial_i (\sigma)(t_1, \dots, t_{n-1})  &= \begin{cases}
		\sigma (t_1, \dots, t_{i-1}, 0, t_i, \dots, t_{n-1}) & i \in \{1, \dots, n\},\\
		\sigma(1- \sum_{j=1}^{n-1} t_j, t_1, \dots, t_{n-1}) \cdot \sigma(v_1)^{-1} & i = 0,
	\end{cases}\\
	s_j (\sigma)(t_1, \dots, t_{n+1}) &= \begin{cases}
		\sigma(t_1, \dots, t_{j-1}, t_j + t_{j+1}, \dots, t_{n+1}) & j \in \{1, \dots, n\},\\
		\sigma(t_2, \dots, t_{n+1}) & j = 0.
	\end{cases}
	\end{aligned}
\end{equation*}
We denote by  $R(A)_\bullet$ the simplicial set $R_\infty(A)_\bullet$. The pointed simplicial sets $R_r (A)_{\bullet}$ and $R (A)_{\bullet}$ are Kan complexes. The fundamental group of $R_r(A)_\bullet$ is given by 
\[\pi_1(R_r(A)_\bullet) = R_r(A)_1/\sim,\]
where $\sim$ denotes the equivalence relation defined by smooth homotopies with fixed end points. The commutator subgroup of $\pi_1(R(A)_\bullet)$ is given by $(F(A)_1/\sim) \subseteq (R(A)_1/\sim)$ the normal closure of the subgroup generated by the elementary matrices 
\begin{equation*}
	e_{ij}(\rho) \in R(A)_\bullet \quad \text{with} \; \rho \in C^\infty([0, 1], M_\infty(A)), \rho(0) = 0.
\end{equation*}
By construction, the relative $K$-groups of $A$ are isomorphic to the homotopy groups of the pointed topological space $|R(A)_\bullet|^+$
	\[K^\mathrm{rel}_n(A) := \pi_n(|R(A)_\bullet|^+), \quad n\ge 1,\]
	where $+$ denotes the Quillen plus construction with respect to the above commutator subgroup of $\pi_1(R(A)_\bullet)$. 

Consider the Banach space
\[\bigwedge\nolimits_n A := A^{\hat \otimes n} / \ker (S).\]
Here the closed subspace $\ker (S) \subseteq A^{\hat \otimes n}$ is the kernel of the continuous map 
\[S: A^{\hat \otimes n} \to A^{\hat \otimes n}, \quad a_1\otimes \cdots \otimes a_n \mapsto \sum_{\tau \in \Sigma_n} \sgn(\tau) a_{\tau^{-1}(1)} \otimes \cdots \otimes a_{\tau^{-1}(n)}.\]
The continuous Lie algebra homology of $A$, denoted by $H^\mathrm{Lie}_\bullet (A)$, is defined to be the homology of the chain complex $(\bigwedge_\bullet A, \delta)$, where $\delta :  \bigwedge_n A \to \bigwedge_{n-1} A$ is the Chevalley-Eilenberg boundary map of $A$ viewed as a Banach Lie algebra
\[a_1 \wedge \cdots \wedge a_n 
\mapsto \sum_{1\le i < j \le n} (-1)^{i+j} [a_i, a_j] \wedge a_1 \wedge \cdots \wedge \hat a_i \wedge \cdots \wedge \hat a_j \wedge \cdots \wedge a_n.\]
Karoubi's relative Chern character $\label{Chernrelative}
	\mathrm{ch}_n^{\mathrm{rel}} : K_n^\mathrm{rel}(A) \to HC_{n-1}(A)$, is a composition
\begin{equation} \label{eq:14n}
	\mathrm{ch}^\mathrm{rel}_n = \mathrm{Tr} \circ \varepsilon \circ L \circ h_n.
\end{equation}
Here the map $h_n$ is the Hurewicz homomorphism associated to the pointed topological space $|R(A)|^+$
\begin{equation*}
	h_n : K_n^\mathrm{rel}(A) = \pi_n(|R(A)_\bullet|^+) \to H_n (|R(A)_\bullet|^+) \cong H_n(R(A)_\bullet).
\end{equation*}
The map $L$ is  the logarithm taking values in the Lie algebra homology 
\begin{equation*}
	\begin{aligned}
		L : H_n(R(A)_\bullet) &\to \lim_{r \to \infty} H_n^\mathrm{Lie} (M_r(A)),
	\end{aligned}
\end{equation*} 
It sends a smooth function $\sigma : \Delta^n \to  \GL_r(A)$ to 
\[\int_{\Delta^n} \dfrac{\partial \sigma}{\partial t_1} \cdot \sigma^{-1} \wedge \cdots \wedge \dfrac{\partial\sigma}{\partial t_n} \cdot \sigma^{-1} dt_1 \cdots dt_n\; .\]
Since $\bigwedge_n A$ is a Banach space, the integral can be defined as a regulated integral for example. The third map is the antisymmetrization
\[\varepsilon : \lim_{r\to \infty} H^\mathrm{Lie}_n(M_r(A)) \to \lim_{r\to \infty} HC_{n-1}(M_r(A)),\]
which is given by the continuous map 
\[
\varepsilon: x_0 \wedge x_1 \wedge \cdots \wedge x_{n-1} \mapsto \frac{(-1)^n}{(n-1)!} \sum_{s\in \Sigma_{n-1}} \sgn(s) x_0 \otimes x_{s(1)} \otimes \cdots \otimes x_{s(n-1)}.
\]
This is the definition of $\varepsilon$ in \cite[\S\,4]{CK}. In \cite[bottom of p. 369]{Kaa11} the normalizing factor $\frac{(-1)^n}{(n-1)!}$ is purposely omitted. 
The last map $\mathrm{Tr}$ is the generalized trace on continuous cyclic homology
\begin{equation}\label{Tr}
	\mathrm{Tr} : \lim_{r \to \infty} HC_{n-1} (M_r(A)) \to HC_{n-1}(A).
\end{equation}
\section{Operator determinants and relative $K$-theory} \label{sec:4}

In this section we review the relation between the de la Harpe, Skandalis approach to determinants for Banach algebras in \cite{HS}, \cite{Har} and Karoubi's theory of the relative Chern character from section \ref{sec:3}. In particular, the Fuglede-Kadison determinant for finite von Neumann algebras $\Nh$ is expressed in terms of $\ch^{\rel}_1$ on $K^{\rel}_1 (\Nh)$. For a unital Banach algebra $A$, recall that $K^{\ttop}_0 (A) = K_0 (A)$ and $K^{\ttop}_i (A) = \pi_{i-1} (\GL (A) , 1)$ for $i \ge 1$. We give more details than in the rest of the paper because we believe that this theory should be better known.

For a group $G$ let $DG = (G,G)$ be the commutator subgroup. A version of the Whitehead lemma shows that $DGL (A) = DGL^0 (A)$ (see e.g. \cite[Lemma 3.1]{Mil71}). Here $\GL^0 (A)$ is the connected component of $1$ in $\GL (A)$. It follows that $DGL (A) \subset \GL^0 (A)$ and hence the map
\[
K_1 (A) = \GL (A) / DGL (A) \twoheadrightarrow K^{\ttop}_1 (A) = \pi_0 (\GL (A)) = \GL (A) / \GL^0 (A)
\]
is surjective. In particular, $K^{\ttop}_1 (A)$ is abelian. Bott periodicity provides natural isomorphisms $K^{\ttop}_i (A) \cong K^{\ttop}_{i+2} (A)$. For $i = 0$ the isomorphism
\begin{equation}
\label{eq:11}
K_0 (A) = K^{\ttop}_0 (A) \iso K^{\ttop}_2 (A) = \pi_1 (\GL^0 (A) , 1) 
\end{equation}
sends the class in $K_0 (A)$ of an idempotent $p \in M_r (A)$ to the homotopy class of the following closed loop with $\alpha_p (0) = 1 = \alpha_p (1)$,
\[
\alpha_p : [0,1] \longrightarrow \GL^0 (A) \; , \; t \mapsto e^{2 \pi i tp} = e^{2 \pi it} p + (1-p) \; 
\]
(see e.g. \cite[sections 1.9-1.11]{Kar78}). Let $\tau$ be a continuous linear map $\tau : A \to \C$ which is tracial, $\tau (ab) = \tau (ba)$ for all $a,b$. It induces a tracial map on $M_{\infty} (A) = \colim_r M_r (A)$ with its inductive limit topology by setting $\tau : M_{\infty} (A) \to \C , \tau ((a_{ij})) = \sum_i \tau (a_{ii})$. Idempotents $p$ and $q$ in $M_{\infty} (A)$ are equivalent if they are conjugate by a matrix in $\GL_{\infty} (A)$. In this case $\tau (p) = \tau (q)$ and hence $\tau$ induces a homomorphism $\tau : K_0 (A) \to \C$.

The de la Harpe, Skandalis determinant for the pair $(A , \tau)$ is obtained as follows. For a piecewise $C^1$-path $\alpha : [t_1 , t_2] \to \GL^0 (A)$ set
\[
\tDelta_{\tau} (\alpha) = \int^{t_2}_{t_1} \tau (\dalpha (t) \alpha (t)^{-1}) \, dt = \tau \int^{t_2}_{t_1} \dalpha (t) \alpha (t)^{-1} \, dt \in \tau (A) \; ,
\]
where $\dalpha(t)$ is the derivative of $\alpha$ with respect to $t$. Since $[t_1 , t_2]$ is compact, $\alpha ([t_1 , t_2]) \subset \GL_r (A) \subset M_r (A)$ for some $r \ge 1$. Hence the integral on the right is defined, for example as a regulated integral of a piecewise continuous function with values in the Banach space $M_r (A)$. According to \cite[Lemma 10]{Har}, under pointwise multiplication of paths $\alpha_0 , \alpha_1 : [t_1 , t_2] \to \GL^0 (A)$ we have
\[
\tDelta_{\tau} (\alpha_0 \alpha_1) = \tDelta_{\tau} (\alpha_0) + \tDelta_{\tau} (\alpha_1) \; .
\]
Moreover, if $H : [t_1 , t_2] \times [0,1] \to \GL^0 (A)$ is a continuous homotopy between $\alpha_0$ and $\alpha_1$ keeping endpoints fixed and with $H ( \_ , s)$ piecewise $C^1$ for each $s$, the map 
\[F(s) := \int_{t_1}^{t_2} \tau \left(\dfrac{\partial H}{\partial t}(t, s) H(t, s)^{-1}\right) dt, \quad s \in [0, 1]\]
is constant. This follows from the tracial property of $\tau$. It then deduces that  $\tDelta_{\tau} (\alpha_0) = \tDelta_{\tau} (\alpha_1)$. Consequently, $\tDelta_{\tau} (\alpha)$ depends only on the homotopy class of $\alpha$ relative to its endpoints.  And it is clear that under concatenation of paths, $\tDelta_{\tau}$ is additive. Finally, for the path $\alpha_p$ above, we have 
\begin{equation}
\label{eq:12}
\tDelta_{\tau} (\alpha_p) = 2 \pi i \tau (p) \; .
\end{equation}
Now, given $x \in \GL^0 (A)$, we choose a piecewise $C^1$-path $\alpha_x : [0,1] \to \GL^0 (A)$ with $\alpha_x (0) = 1 , \alpha_x (1) = x$ and set
\[
\Delta_{\tau} (x) \equiv \tDelta_{\tau} (\alpha_x) \equiv \int^1_0 \tau (\dalpha_x (t) \alpha_x (t)^{-1}) \, dt \, \mod 2 \pi i \tau (K_0 (A)) \; .
\]
This is well defined because if $\alpha_x , \alpha'_x : [0,1] \to \GL^0 (A)$ are paths from $1$ to $x$, then $\alpha = \alpha^-_x \verk \alpha'_x$ is a closed path and hence defines an element of $\pi_1 (\GL^0 (A) , 1)$, where $\alpha_x^{-}$ is the same path as $\alpha_x$ but with inverse direction. Thus
\[
\tDelta_{\tau} (\alpha'_x) - \tDelta_{\tau} (\alpha_x) = \tDelta_{\tau} (\alpha) \in \tDelta_{\tau} (\pi_1 (\GL^0 (A) , 1)) = 2 \pi i \tau (K_0 (A)) \; .
\]
Here, the last equality follows from \eqref{eq:11} and \eqref{eq:12}. Noting that $\tDelta_{\tau} (e^{ta}) = \tau (a)$ where $a \in A$ and $t \mapsto e^{ta}$ is viewed as a path on $[0,1]$, we get in conclusion \cite[Theorem 12]{Har}.

\begin{theorem}
\label{t2.1}
For $(A, \tau)$ as above, the map
\[
\Delta_{\tau} : \GL^0 (A) \longrightarrow \C / 2 \pi i \tau (K_0 (A))
\]
is a homomorphism of groups with image $\tau (A) / 2 \pi i \tau (K_0 (A))$. 
\end{theorem}

The Fuglede-Kadison determinant which we recalled in the introduction is a special case of $\Delta_{\tau}$. Namely, a von Neumann algebra $\Nh$ is a Banach algebra and it is known that $\GL (\Nh) = \GL^0 (\Nh)$ is (path-)connected. (Use the polar decomposition and the \textit{measurable} functional calculus in $\Nh$ which guarantees the existence of a logarithm $C$ of a unitary operator $U$ to connect $1$ to $U$ via $\exp tC$ for $t \in [0,1]$.) Hence $K^{\ttop}_1 (\Nh) = 0$ and by Bott periodicity $K^{\ttop}_{2i + 1} (\Nh) = 0$ for $i \ge 0$. Assume that $\Nh$ admits a finite positive trace $\tau : \Nh \to \C$. Then we have $\tau (a^*) = \overline{\tau (a)}$ for $a \in \Nh$. Since every class in $K_0 (\Nh)$ is represented by a self-adjoint projection, it follows that $\tau (K_0 (\Nh)) \subset \R$. The Fuglede-Kadison determinants for the von Neumann algebras $M_r (\Nh)$ with trace $\tau ((a_{ij})) = \sum_i^r \tau (a_{ii})$ assemble into a homomorphism $\det_{\tau} : \GL (\Nh) = \colim_r \GL_r(\Nh) \to \R^{> 0}$. The resulting homomorphism $\det_{\tau} : K_1 (\Nh) \to \R^{> 0}$ is related to $\Delta_{\tau}$ as follows, \cite[Corollary 14]{Har}.

\begin{theorem}
\label{t2.2}
For a von Neumann algebra $\Nh$ with a finite positive trace $\tau$ there is a commutative diagram
\[
\xymatrix{
\GL (\Nh) \ar@{->>}[d] \ar@{=}[r] & \GL^0 (\Nh) \ar[r]^-{\Delta_{\tau}} & \C / 2 \pi i \tau (K_0 (\Nh)) \ar[d]^{\RRe} \\
K_1 (\Nh) \ar[rr]^{\log \det_{\tau}} & & \R \; .
}
\]
\end{theorem}

The interpretation of the de la Harpe, Skandalis determinant $\Delta_{\tau}$ of Theorem \ref{t2.1} in terms of relative $K$-theory is as follows, \cite[5.1]{CK}, \cite[3.0]{Kaa11}, \cite[\S\,5]{HKS}. Recall that for a Fr\'echet algebra $A$, the $0$-th topological cyclic homology is given by $HC_0 (A) = A / \Imm b$ where the continuous linear differential $b : A \hat{\otimes}_{\pi} A \to A$ is uniquely determined by the formula $b (a_1 \otimes a_2) = [a_1 , a_2]$. Thus we have\footnote{In the pioneering works \cite{Kar87} and \cite[5.1]{CK} contrary to the improved theory of \cite{Kar97}, the maximal Hausdorff quotients of the homologies are used. For example, in these references $HC_0 (A)$ denotes the quotient of $A$ by the closure of $\Imm b$ which is the same as the quotient of $A$ by the closure of $[A,A]$.} $[A , A] \subset \Imm b \subset \overline{[A,A]}$. Consider the diagram \eqref{diag1} in low degrees where $w = \ch^{\rel}_1 \verk \delta$
\begin{equation} \label{eq:13}
\vcenter{\xymatrix{
K^{\ttop}_2 (A) \ar[d]^{\ch^{\ttop}_2} \ar[r]^{\delta} \ar[dr]^w & K^{\rel}_1 (A) \ar[d]^{\ch^{\rel}_1} \ar[r]^{\theta} & K_1 (A) \ar[d]^{\ch^-_1} \ar[r]^v & K^{\ttop}_1 (A) \ar[d] \\
HP_2 (A) \ar[r]^S & HC_0 (A) \ar[r]^{\partial} & HC^-_1 (A) \ar[r] & HP_1 (A) \; .
}}
\end{equation}
%
A tracial functional $\tau : A \to \C$ defines a linear map $\tau : HC_0 (A) \to \C$. As a special case of diagram \eqref{eq:27} below, we have $\tau (w (K^{\ttop}_2 (A))) = 2 \pi i \tau (K_0 (A))$. Here $\tau (K_0 (A))$ as defined above is really $\tau (ch^{\ttop}_0 (K_0 (A)))$. 

The map $v$ can be identified with the natural projection 
\[
\GL (A) / DGL (A) = \GL (A) / DGL^0 (A) \to \GL_ (A) / \GL^0 (A) \; .
\]
Hence its kernel is
\[
\Ker v = \GL^0 (A) / DGL^0 (A) \; .
\]
Since $\Delta_{\tau}$ in Theorem \ref{t2.1} is a homomorphism, it may be viewed as a map
\[
\Delta_{\tau} : \Ker v \longrightarrow \C / 2 \pi i \tau (K_0 (A)) \; .
\]

\begin{theodef}
\label{t2.3}
The following diagram commutes for any Banach algebra $A$ and any tracial funtional $\tau : A \to \C$,
\begin{equation} \label{eq:14}
\vcenter{\xymatrix{
K^{\rel}_1 (A) / \delta (K^{\ttop}_2 (A)) \ar[d]^{\ch^{\rel}_1} \ar[r]^{\overset{\theta}{\simeq}} & \Ker v \subset K_1 (A) \ar[d]^{\Delta_{\tau}} \\
HC_0 (A) / w (K^{\ttop}_2 (A)) \ar[r]^{\tau} & \C / 2 \pi i \tau (K_0 (A)) \; .
}}
\end{equation}
For Fr\'echet algebras $A$ with a continuous trace $\tau$, we define $\Delta_{\tau}$ by diagram \eqref{eq:14}.
\end{theodef}

\begin{proof} Let $A$ be a Banach algebra. Recall the explicit description of $\ch^{\rel}_1$ at the end of section \ref{sec:3}. We have:
\[
K^{\rel}_1 (A) = \colim_r (R_r (A)_1 / \sim) / \pi (F_r (A)_1) \; ,
\]
where
\[
R_r (A)_1 = \{ \sigma \in C^{\infty} ([0,1], \GL_r (A)) \mid \sigma (0) = 1 \} \quad \text{and} \quad F_r (A)_1 = (R_r (A)_1 , R_r (A)_1) \; .
\]
Here, $\sim$ is the equivalence relation on $R_r (A)_1$ given by smooth homotopies with fixed endpoints and $\pi : R_r (A)_1 \to R_r (A)_1 / \sim$ is the projection. According to formula \eqref{eq:14n} the relative Chern character $\ch^{\rel}_1 : K^{\rel}_1 (A) \to HC_0 (A)$ is induced by the following maps for $r \ge 1$ where $\Tr$ is the ordinary matrix trace
\begin{equation}
	\label{eq:15}
	R_r (A)_1 \longrightarrow A \; , \; \sigma \longmapsto -\Tr \int^1_0 \dsigma (t) \sigma (t)^{-1} \, dt \; .
\end{equation} 
Hence we may write
\[
\ch_1^\mathrm{rel}([\sigma]) = - \left[\mathrm{Tr} \int_0^1 \dot{\sigma}(t) \sigma(t)^{-1} dt\right],
\]
The map $\theta$ in diagram \eqref{eq:13} sends the class $[\sigma]$ of $\sigma \in R_r (A)_1$ to $\theta ([\sigma)] = \sigma (v_1)^{-1} = \sigma (1)^{-1}$. We therefore have
\begin{align*}
\Delta_{\tau} (\theta ([\sigma])) & = -\Delta_{\tau} (\sigma (1)) = -\tDelta_{\tau} (\sigma) = -\tau \int^1_0 \dsigma (t) \sigma (t)^{-1} \, dt \\
& := -\tau \Big( \Tr \int^1_0 \dsigma (t) \sigma (t)^{-1} \, dt\Big) = \tau (\ch^{\rel}_1 ([\sigma])) \; .
\end{align*}
Hence diagram \eqref{eq:14} commutes. To understand the working of Karoubi's theory in the simplest case, it is instructive to see directly that for homotopic $\sigma_0 , \sigma_1 \in R_p (A)_1$ the images under the map \eqref{eq:15} differ by an element in $\Imm b$. We follow the proof of \cite[Proposition 4.5]{HKS}. Let $H : [0,1] \times [0,1] \to \GL_r (A)$ be a smooth homotopy between $\sigma_0$ and $\sigma_1$ with fixed endpoints. Thus $H (t,0) = \sigma_0 (t) , H (t,1) = \sigma_1 (t)$ for all $t$ and $H (0,s) , H (1,s)$ are constant in $s$. Noting that $\frac{\partial H}{\partial s} (1,s) = 0 = \frac{\partial H}{\partial s} (0,s)$ and $\frac{\partial H^{-1}}{\partial t} = -H^{-1} \frac{\partial H}{\partial t} H^{-1}$ we get the following equations in $M_r (A)$:
\begin{eqnarray*}
\lefteqn{
\int^1_0 \dsigma_1 (t) \sigma_1 (t)^{-1} \, dt - \int^1_0 \dsigma_0 (t) \sigma_0 (t)^{-1} \, dt}\\
& = & \int^1_0 \Big( \frac{\partial H}{\partial t} (t,1) H (t,1)^{-1} - \frac{\partial H}{\partial t} (t,0) H (t,0)^{-1} \Big) \, dt\\
& & - \int^1_0 \Big( \frac{\partial H}{\partial s} (1,s) H (1,s)^{-1} - \frac{\partial H}{\partial s} (0,s) H (0,s)^{-1} \Big) \, ds \\
& = & \int^1_0 \int^1_0 \frac{\partial}{\partial s} \Big( \frac{\partial H}{\partial t} H^{-1} \Big) \, ds \, dt - \int^1_0 \int^1_0 \frac{\partial}{\partial t} \Big( \frac{\partial H}{\partial s} H^{-1} \Big) \, dt \, ds \\
& = & \int^1_0 \int^1_0 \Big[ \frac{\partial H}{\partial s} H^{-1} , \frac{\partial H}{\partial t} H^{-1} \Big] \, dt \, ds \; .
\end{eqnarray*}
Since
\[
\Tr [a,b] = \sum_{i,j} [a_{ij} , b_{ji}] \quad \text{for matrices} \; a = (a_{ij}) , b = (b_{ij}) \; ,
\]
it follows that
\begin{eqnarray*}
\lefteqn{\Tr \int^1_0 \dsigma_1 (t) \sigma_1 (t)^{-1} \, dt - \Tr \int^1_0 \dsigma_0 (t) \sigma_0 (t)^{-1} \, dt} \\
& = & \sum_{i,j} \int^1_0 \int^1_0 \Big[ \Big( \frac{\partial H}{\partial s} H^{-1} \Big)_{ij} , \Big( \frac{\partial H}{\partial t} H^{-1} \Big)_{ji} \Big] \, dt \, ds\\
& = & b (l (H)) \; ,
\end{eqnarray*}
where
\[
l (H) = \sum_{i,j} \int^1_0 \int^1_0 \Big( \frac{\partial H}{\partial s} H^{-1} \Big)_{ij} \otimes \Big( \frac{\partial H}{\partial t} H^{-1} \Big)_{ji} \, dt \, ds \in A \hat{\otimes}_{\pi} A \; .
\]
Here the integral over $[0,1] \times [0,1]$ is a regulated integral with values in the Banach space $A \hat{\otimes}_{\pi} A$. It follows that the relative Chern character on $K^{\rel}_1 (A)$ is well defined as a homomorphism into $HC_0 (A) = A / \Imm b$ and then of course into $\overline{HC_0} (A) = A / \overline{\Imm b} = A / \overline{[A,A]}$ as in the first works \cite{Kar87} and \cite{CK}. However, it does not exist as a map into algebraic cyclic homology $HC_0 (A^{\delta}) = A / [A,A]$ of $A$ with the discrete topology. The relative Chern character is inherently topological in nature.
\end{proof}

\begin{cor} \label{t2.4}
1) For a finite von Neumann algebra $\Nh$ with a finite positive trace $\tau : \Nh \to \C$, we have a commutative diagram
\begin{equation}
\label{eq:16}
\vcenter{\xymatrix{
K^{\rel}_1 (\Nh) / \delta K^{\ttop}_2 (\Nh) \ar[d]_{\ch^{\rel}_1} \ar@{->>}[r]^-{\theta} & K_1 (\Nh) \ar[dd]^{\log \ddet_{\tau}} \\
HC_0 (\Nh) / w K^{\ttop}_2 (\Nh) \ar[d]_{\tau} & & \\
\C / 2 \pi i \tau (K_0 (\Nh)) \ar[r]^-{\RRe} & \R \; .
}}
\end{equation}
\end{cor}

\begin{proof} 
This follows from Theorems \ref{t2.2} and \ref{t2.3}, noting that $\tau (w K^{\ttop}_2 (\Nh)) = \\
2 \pi i \tau (K_0 (\Nh)) \subset 2 \pi i \R$. Here as in section \ref{sec:2}, we view the trace $\tau$ as a homomorphism on $K_0 (A)$ by setting $\tau ([p]) = \tau (\Tr \, p)$ where $p^2 = p$ is a projector on $A^n$ for some $n$. 
\end{proof}

In particular, for $\Nh = L^{\infty} (\Omega , \mu)$, where $\Omega$ is a standard measure space with probability measure $\mu$ and for the trace $\tau (f) = \int_{\Omega} f \, d\mu$, we find for $a \in K_1 (\Nh)$, that
\[
\log \ddet_{\tau} (a) = \RRe \int_{\Omega} \ch^{\rel}_1 (\theta^{-1} (a)) \, d \mu \; .
\]
This gives an interpretation for the integral
\[
\int_{\Omega} \log |\det (a)|^{1 / r} d \mu = \log \ddet_{\tau} (a) \quad \text{where} \; a \in \GL_r (L^{\infty} (\Omega , \mu)) 
\]
in terms of $\ch^{\rel}_1$. In particular, for $M = T^d$ and the Haar measure $\mu$ of $T^d$, the (logarithmic) Mahler measure $m (a)$ of $a \in L^{\infty} (T^d , \mu)^{\times}$ can be written as follows
\begin{equation}
\label{eq:17}
m (a) = \int_{T^d} \log |a| d\mu = \log \ddet_{\tau} (a) = \RRe \int_{T^d} \ch^{\rel}_1 (\theta^{-1} (a)) \, d \mu \; .
\end{equation}

\section{Passing from $\ch^{\rel}_1$ to $\ch^{\rel}_{d+1}$} \label{sec:5}

In section \ref{sec:4}, we saw that the Fuglede-Kadison determinant, and in particular the Mahler measure can be expressed in terms of Karoubi's relative Chern character map $\ch^{\rel}_1 : K^{\rel}_1 (A) \to HC_0 (A)$ evaluated on a cyclic cocyle, a trace $\tau \in HC^0 (A)$. In this section we explain for commutative Fr\'echet algebras $A$ how to pass from $K^{\rel}_1 (A)$ to $K^{\rel}_{d+1} (A)$ for $d \ge 1$ by using a module structure under algebraic $K$-theory. The relevant result is the following consequence of \cite[Th\'eor\`eme 6.28]{Kar87}.

\begin{theorem}[Karoubi] \label{t5.1}
Let $A$ be a commutative Fr\'echet algebra. The graded abelian group $K^{\rel}_* (A)$ is a module under the graded commutative ring $K_* (A)$. The cyclic homology $HC_* (A)$ is a graded module under the graded commutative ring $HC^-_* (A)$. The maps $K^{\rel}_* (A) \to K_* (A)$ and $HC_{*-1} (A) \to HC^-_* (A)$ are $K_* (A)$-- and $HC^-_* (A)$ equivariant. The Chern characters $\ch^-_i : K_i (A) \to HC^-_i (A)$ and $\ch^{\rel}_j : K^{\rel}_j (A) \to HC_{j-1} (A)$ satisfy the following formula for $x \in K_i (A)$ and $y \in K^{\rel}_j (A)$
\begin{equation}
\label{eq:24}
\ch^{\rel}_{i+j} (xy) = \ch^-_i (x) \ch^{\rel}_j (y) \quad \text{in} \; HC_{i+j-1} (A) \; .
\end{equation}
\end{theorem}

\textbf{Comments on Theorem \ref{t5.1}} Strictly speaking, Karoubi proved the result only for Banach algebras, but in \cite[6.19]{Kar87} he remarks that the adaptation to the simplicial case (for Fr\'echet algebras) poses no problems. \\
The target groups of $\ch^-_i$ and $\ch^{\rel}_j$ in Karoubi's notation are $\Hh^*_A (S^i)$ and $\Hh^{'*}_A (S^j)$. Using the formula for $\Hh^{2n} (X,A)$ in \cite[p. 77]{Kar87} for $X = S^i$ one sees that $\Hh^*_A (S^i)$ is isomorphic to negative cyclic homology and that the map $\ch^{\rel}$ into $\Hh^{'*}_A (S^j)$ maps into cyclic homology which is naturally a subgroup of $\Hh^{'*}_A (S^j)$. More precisely, Karoubi obtains the maximal Hausdorff quotients of these groups, but this is not a problem for $A = C^{\infty} (M)$, where $M$ is a compact manifold. In this case the cyclic and negative cyclic homologies are Hausdorff as one sees from the explicit formulas in section \ref{sec:3} noting that the image of the de~Rham differential is closed in the Fr\'echet topology of the space of smooth forms on $M$ because the de~Rham complex is elliptic. Karoubi's description of algebraic $K$-theory is more explicit than Quillen's and correspondingly for $K^{\rel}$. In \cite[Theorem 3.11 and Proposition 3.12]{Kar87} he shows that his $K$-groups together with their product agree with Quillen's $K$-groups equipped with the Loday product. He defines the $K_*$-module structure of $K^{\rel}_*$ at the beginning of the proof of \cite[Th\'eor\`eme 6.28]{Kar87}. It would be a useful exercise to reprove Theorem \ref{t5.1} for general Fr\'echet algebras using a Loday style product to define the $K_*$-module structure on $K^{\rel}_*$ and the $HC^-_*$-module structure on $HC_*$ extended from the algebraic case (\cite[5.1.14, 5.1.15]{Loday}, \cite[\S\,2]{BO}) to topological algebras. The method in \cite{Gi} used to prove that $\ch^-_*$ is a ring homomorphism will certainly be helpful.

Generalizing Theorem / Definition \ref{t2.3}, we define a higher de la Harpe Skandalis determinant for Fr\'echet algebras $A$ using Karoubi's commutative diagram 
\begin{equation}
\label{eq:25}
\vcenter{
\xymatrix{
K^{\ttop}_{d+2} (A) \ar[r]^{\delta} \ar[d]_{\ch^{\ttop}_{d+2}} \ar[dr]^w & K^{\rel}_{d+1} (A) \ar[r]^{\theta} \ar[d]_{\ch^{\rel}_{d+1}} & K_{d+1} (A) \ar[r]^v \ar[d]_{\ch^-_{d+1}} & K^{\ttop}_{d+1} (A) \ar[d]_{\ch^{\ttop}_{d+1}} \\
HP_{d+2} (A) \ar[r]^S & HC_d (A) \ar[r]^{\partial} & HC^-_{d+1} (A) \ar[r] & HP_{d+1} (A)
}}
\end{equation}
The determinant depends on a cyclic cohomology class $\tau^d \in HC^d (A)$ which we view as a linear form on $HC_d (A)$ and via composition $\tau^d : HP_d (A) \xrightarrow{\lambda} HC_d (A) \xrightarrow{\tau^d} \C$ also on $HP_d (A)$. 

\begin{definition}
For a Fr\'echet algebra $A$ and a cyclic cocycle $\tau^d \in HC^d (A)$ the $d$-th de la Harpe Skandalis determinant $\Delta_{\tau^d}$ is defined by the commutative diagram
\begin{equation}
\label{eq:26}
\vcenter{
\xymatrix{
K^{\rel}_{d+1} (A) / \delta (K^{\ttop}_{d+2} (A)) \ar[r]^{\overset{\theta}{\sim}} \ar[d]_{\ch^{\rel}_{d+1}} & \Ker (K_{d+1} (A) \to K^{\ttop}_{d+1} (A)) \ar[d]^{\Delta_{\tau^d}} \\
HC_d (A) / 2 \pi i \, \lambda (\ch^{\ttop}_d (K^{\ttop}_d (A))) \ar[r]^{\tau^d} & \C / 2 \pi i \tau^d (\ch^{\ttop}_d (K^{\ttop}_d (A))) \; .
}}
\end{equation}
\end{definition}

A priori, the indeterminacy of $\Delta_{\tau^d}$ lies in $\tau^d (w (K^{\ttop}_{d+2} (A))$, but note that $w (K^{\ttop}_{d+2} (A)) = 2 \pi i \lambda (\ch^{\ttop}_d (K^{\ttop}_d (A))$ and hence $\tau^d (w (K^{\ttop}_{d+2} (A)) = 2 \pi i \tau^d (\ch^{\ttop}_d (K^{\ttop}_d (A)))$ because of the following commutative diagram \cite{C}
\begin{equation}
\label{eq:27}
\vcenter{
\xymatrix{
K^{\ttop}_d (A) \ar[r]^{\ch^{\ttop}_d} \ar[d]^{\wr}_{\text{Bott}} & HP_d (A) \ar[r]^{\lambda} \ar[d]^{2 \pi i S^{-1}} & HC_d (A) \ar[d]^{2 \pi i} \\
K^{\ttop}_{d+2} (A) \ar[r]^{\ch^{\ttop}_{d+2}} & HP_{d+2} (A) \ar[r]^S & HC_d (A)
}}
\end{equation}
The following consequence of Theorem \ref{t5.1} allows us to express the original de la Harpe Skandalis determinant for certain traces by higher determinants if $A$ is commutative.

\begin{cor} \label{t5.2}
For a commutative Fr\'echet algebra $A$ and $x \in K_d (A) , \tau^d \in HC^d (A)$, consider the trace $\tau^0_x \in HC^0 (A)$ defined by the following linear form where we use the $HC^-_* (A)$-module structure on $HC_* (A)$: 
\[
\tau^0_x : A = HC_0 (A) \xrightarrow{\ch^- (x) \scriptstyle{\bullet}} HC_d (A) \xrightarrow{\tau^d} \C \; .
\]
Then $2 \pi i \tau^0_x (K_0 (A))$, the indeterminacy of $\Delta_{\tau^0_x}$ is contained in $2 \pi i \tau^d (\ch^{\ttop}_d (K^{\ttop}_d (A))$, the indeterminacy of $\Delta_{\tau^d}$. Moreover, for $f \in \Ker (K_1 (A) \xrightarrow{v} K^{\ttop}_1 (A))$ there is the following formula
\begin{equation}
\label{eq:28n}
\Delta_{\tau^0_x} (f) = \Delta_{\tau^d} (xf) \quad \text{in} \; \C / 2 \pi i \, \tau^d (\ch^{\ttop}_d (K^{\ttop}_d (A))) \; .
\end{equation}
\end{cor}

\begin{proof}
Let $y \in K^{\rel}_1 (A)$. Using Theorem \ref{t5.1} we get
\begin{align*}
\tau^0_x (\ch^{\rel}_1 (y)) & = \tau^d (\ch^-_d (x) \ch^{\rel}_1 (y)) \\
& = \tau^d (\ch^{\rel}_{d+1} (xy)) \; .
\end{align*}
In particular, for $y = \delta (\Bott \cdot \varepsilon) \in \delta (K^{\ttop}_2 (A)) \subset K^{\rel}_1 (A)$, where $\Bott \cdot \varepsilon$ denotes the element of $K^{\ttop}_2 (A)$ obtained by applying the Bott periodicity isomorphism to $\varepsilon \in K_0 (A)$, we find
\begin{align*}
2 \pi i \tau^0_x (\varepsilon) & := 2 \pi i \tau^0_x (\ch^{\ttop}_0 (\varepsilon)) \overset{\eqref{eq:27}}{=} \tau^0_x (S (\ch^{\ttop}_2 (\Bott \cdot \varepsilon))) \\
& \overset{\eqref{eq:25}}{=} \tau^0_x (\ch^{\rel}_1 (\delta (\Bott \cdot \varepsilon)) \\
& = \tau^d (\ch^{\rel}_{d+1} (x \cdot \delta (\Bott \cdot \varepsilon)) \; .
\end{align*}
Since $\theta (x \delta (\Bott \cdot \varepsilon)) = x \theta (\delta (\Bott \cdot\varepsilon)) = 0$ by $K_* (A)$-equivariance of the map $\theta$, it follows that
\[
x\delta (\Bott \cdot \varepsilon) \in \delta (K^{\ttop}_{d+2} (A)) \; ,
\]
and therefore
\[
2 \pi  i \tau^0_x (\varepsilon) \in \tau^d (w (K^{\ttop}_{d+2} (A)) = 2 \pi i \tau^d (\ch^{\ttop}_d (K^{\ttop}_d (A))) \; .
\]
For the indeterminacies of $\Delta_{\tau^0_x}$ and $\Delta_{\tau^d}$ we therefore have
\[
2 \pi i \tau^0_x (K_0 (A)) \subset 2 \pi i \tau^d (\ch^{\ttop}_d (K^{\ttop}_d (A)) \; .
\]
For $f \in \Ker (K_1 (A) \xrightarrow{v} K^{\ttop}_1 (A))$, choose $y \in K^{\rel}_1 (A)$ with $f = \theta (y)$. Again using that $\theta$ is $K_* (A)$-equivariant, we find that
\[
xf = x \theta (y) = \theta (xy) \in \Ker (K_{d+1} (A) \xrightarrow{v} K^{\ttop}_{d+1} (A)) \; .
\]
Applying the definitions of $\Delta_{\tau^0_x}$ and $\Delta_{\tau^d}$ and the formula from the beginning of the proof, we get
\begin{align*}
\Delta_{\tau^0_x} (f) & = \tau^0_x (\ch^{\rel}_1 (y)) = \tau^d (\ch^{\rel}_{d+1} (xy)) \\
& = \Delta_{\tau^d} (\theta (xy)) = \Delta_{\tau^d} (xf) \; .
\end{align*}
\end{proof}

We now apply Corollary \ref{t5.2} to the Fr\'echet algebra $A = C^{\infty} (T^d)$. For the element $x \in K_d (A)$ we take the image under the assembly map of the fundamental class of $T^d$ with the orientation given by the ordering $t_1 , \ldots , t_d$, see \cite[12.4.1]{Loday}
\[
\xymatrix{
H_d (T^d , \Q) \ar@{=}[r] & H_d (\Z^d , \Q) \ar[r]^{\text{assembly}} & K_d (\Z [\Z^d])_{\Q} \ar[r] & K_d (C^{\infty} (T^d))_{\Q} \ar[d]_{ch^-_d} \\
[T^d] \ar@{|->}[r] & \frac{dt_1}{t_1} \wedge \ldots \wedge \frac{dt_d}{t_d} \ar@{}[r]|{\in} & \Ah^d (T^d) \ar@{=}[r] & HC^-_d (C^{\infty} (T^d))
}
\]
Here $\Ah^d (T^d)$ is the space of smooth $d$-forms on $T^d$ (they are closed). The fact that $[T^d]$ goes to $\frac{dt_1}{t_1} \wedge \ldots \wedge \frac{dt_d}{t_d}$ follows from multiplicativity of all maps and the simple case $d = 1$. For the module structure of cyclic homology under negative cyclic homology of commutative Fr\'echet algebras the pairing
\[
\xymatrix{
HC_0 (C^{\infty} (T^d)) \times HC^-_d (C^{\infty} (T^d)) \ar[r] \ar@{=}[d] & HC_d (C^{\infty} (T^d)) \ar@{=}[d] \\
C^{\infty} (T^d) \times \Omega^d (T^d) \ar[r] & \Ah^d (T^d) / d \Ah^{d-1} (T^d) \oplus H^{d-2} (T^d) \oplus \ldots
}
\]
is given by sending $(f, \omega)$ to the class of $f\omega$ in $H_d (T^d)$. Consider the cyclic cocycle $\tau^d \in HC^d (T^d)$, given by 
\[
\tau^d (c) = (2 \pi i)^{-d} \int_{T^d} c_d \quad \text{for} \; c \in HC_d (T^d) \; ,
\]
where $c_d$ is the component of $c$ in $\Ah^d (T^d) / d \Ah^{d-1} (T^d)$. For $x = \{ T^d \}$ defined as the image of $[T^d]$ under the assembly map into $K_d (C^{\infty} (T^d))$ we obtain with $\alpha \in C^{\infty} (T^d)$
\begin{align*}
\tau^0_x (\alpha) & = \tau^d (ch_d (x) \alpha) = \tau^d \Big( \alpha \frac{dt_1}{t_1} \wedge \ldots \wedge \frac{dt_d}{t_d} \Big) = (2 \pi i)^{-d} \int_{T^d} \alpha \frac{dt_1}{t_1} \wedge \ldots \wedge \frac{dt_d}{t_d} \\
& = \int_{T^d} \alpha \, d\mu = \tau \; , \; \text{the standard positive trace on} \; C^{\infty} (T^d) \subset \Nh \Z^d = L^{\infty} (T^d , \mu) \; .
\end{align*}

The projection of $\ch^{\ttop}_d (K^{\ttop}_d (C^{\infty} (T^d)))$ to the summand $H^d (T^d , \C)$ of $HP_d (C^{\infty} (T^d))$ is equal to the image of the topological Chern character $K^{-d}_{\ttop} (T^d) \to H^d (T^d , \Q (d))$  where $\Q (d) = (2 \pi i)^d \Q$. On $H^d (T^d , \Q (d))$, our trace $\tau^d = (2 \pi i)^{-d} \int_{T^d}$ takes values in $\Q$. Using Corollary \ref{t5.2} we therefore obtain the following result.

\begin{cor} \label{t5.3}
We have
\[
2 \pi i \tau^d (\ch^{\ttop}_d (K^{\ttop}_d (C^{\infty} (T^d))) \subset 2 \pi i \Q \; .
\]
Moreover, for $f \in \Ker (K_1 (C^{\infty} (T^d)) \to K^{\ttop}_1 (C^{\infty} (T^d))) \otimes \Q$ we have
\[
\Delta_{\tau} (f) = \Delta_{\tau^d} ( \{ T^d \} \cdot f) \quad \text{in} \; \C / 2 \pi i \Q \; .
\]
\end{cor}

Using functoriality of Karoubi's diagram \eqref{eq:13} under the continuous inclusion $C^{\infty} (T^d) \subset L^{\infty} (T^d , \mu)$, formula \eqref{eq:17} and Corollary \ref{t5.3} we get: 

\begin{cor} \label{t5.4} 
For $f \in \Ker (K_1 (C^{\infty} (T^d)) \to K^{\ttop}_1 (C^{\infty} (T^d))) \otimes \Q$ we have
\begin{equation}
\label{eq:29n}
m (f) = \log \ddet_{\tau} (f) = \RRe \Delta_{\tau} (f) = \RRe \Delta_{\tau^d} ( \{ T^d \} \cdot f) \; .
\end{equation}
\end{cor}

\begin{rem}
An element $f \in C^{\infty} (T^d)^{\times}$ gives rise to an element in the kernel of the map $K_1 (C^{\infty} (T^d))_{\Q} \to K^{\ttop}_1 (C^{\infty} (T^d))_{\Q}$ if and only if the class of $(df) f^{-1}$ in $H^1 (T^d , \C)$ is zero. This can always be achieved by multiplying $f$ with $t^{\nu_1}_1 \cdots t^{\nu_d}_d f$ for a uniquely determined $\nu \in \Z^d$. It is equivalent to say that for $f$ a smooth logarithm $\log f$ exists. 
\end{rem}

Because of formula \eqref{eq:29n} we now have a conceptual understanding why for $f \in \Z [t^{\pm 1}_1 , \ldots , t^{\pm 1}_d]$ with $f (z) \neq 0$ for all $z \in T^d$, the logarithmic Fuglede-Kadison-determinant of $f$ is related to higher algebraic $K$-theory and a higher Chern-character, namely $ch^{\rel}_{d+1}$. The relation with the Deligne-Beilinson regulator on algebraic $K$-theory given by formula \eqref{eq:2} can now be established conceptually by using results in Tamme's thesis which compare Karoubi's and Beilinson's higher Chern class maps. This will be done in the next section.
\section{Passing from $ch^{\rel}_{d+1}$ to $r_{\Dh}$, the work of Tamme} \label{sec:6}
We use Tamme's comparison theorems in \cite{T1} and \cite{Tam12} to show that formula \eqref{eq:2} follows from formula \eqref{eq:29n} via natural commutative diagrams. First some background. Let $X = \spec R$ be a smooth affine scheme of finite type over $\C$. By definition $K_i (X) = K_i (R)$ and $K^{-i}_{\ttop} (X) = \pi_i (BU^X)$. Here $BU^X$ is the classifying space of $U^X$, the topological group of continuous maps from the complex manifold $X (\C)$ to the infinite unitary group with the compact open topology. By \cite[3.2]{T1}, there is a natural map $|B_{\bullet} \GL (R)|^+ \to BU^X$. Tamme's relative $K$-groups $K^{\rel}_i (X)$ are defined as the $i$-th homotopy group of the homotopy fibre of this map. Hence there is a natural exact sequence
\[
\ldots \longrightarrow K_{i+1} (X) \longrightarrow K^{-i-1}_{\ttop} (X) \longrightarrow K^{\rel}_i (X) \longrightarrow K_i (X) \longrightarrow K^{-i}_{\ttop} (X) \longrightarrow \ldots 
\]
Let $F^{\bullet}$ be the Hodge filtration on $H^q (X (\C) , \C)$ and let $H^q_{\Dh} (X , \R (n))$ be the (algebraic) real Deligne cohomology grups of $X / \C$. They fit into an exact sequence \cite[Corollary 2.10, b)]{EV},
\begin{equation}
\label{eq:28}
\ldots \to H^q (X (\C), \C) / F^n \xrightarrow{\alpha} H^{q+1}_{\Dh} (X , \R (n)) \to H^{q+1} (X (\C) , \R (n)) \to \ldots
\end{equation}
In \cite[3.3]{T1}, Tamme defined relative Chern characters $ch^{\rel}_{n,i}$ in the spirit of Karoubi's work and proved that they fit into a commutative diagram\footnote{The normalization of $ch^{\rel}_{n,i}$ in \cite{T1} is compatible with Karoubi's. In \cite{Tam12} Tamme used a different normalization so that diagram  (31) in \cite{Tam12} does not have the sign $(-1)^{n-1}$.} \cite[Theorem 3.11]{T1}
\begin{equation} \label{eq:29}
\vcenter{
\xymatrix{
K^{\rel}_i (X) \ar[r]^{\theta} \ar[d]_{(-1)^{n-1} ch^{\rel}_{n,i}} & K_i (X) \ar[r] \ar[d]^{ch^{\Dh}_{n,i}} & K^{\ttop}_{-i} (X) \ar[d]^{ch^{\ttop}_{n,i}} \\
H^{2n-i-1} (X (\C) , \C) / F^n \ar[r]^{\alpha} & H^{2n-i}_{\Dh} (X , \R (n)) \ar[r] & H^{2n-i} (X (\C) , \R (n))
}}
\end{equation}
On the other hand, at the end of section 3.3 in \cite{T1}, Tamme shows that for a version of Karoubi's relative Chern character, the following diagram commutes:
\begin{equation}
\label{eq:30}
\vcenter{
\xymatrix{
K^{\rel}_i (X) \ar[r] \ar[d]_{ch^{\rel}_{n,i}} & K^{\rel}_i (C^{\infty}  (X)) \ar[d]^{ch^{\rel}_{n,i}} \\
H^{2n-i-1} (X (\C) , \C) / F^n \ar[r]^{\beta} & H^{2n-i-1} (\Gamma (X (\C) , \Ah^{\bullet < n})) \; .
}}
\end{equation}
Let $M$ be a smooth compact submanifold of $X (\C)$. The same construction as the one in \cite{T1} for diagram \eqref{eq:30} gives a commutative diagram:
\begin{equation}
\label{eq:31}
\vcenter{
\xymatrix{
K^{\rel}_i (X) \ar[r] \ar[d]_{ch^{\rel}_{n,i}} & K^{\rel}_i (C^{\infty} (M)) \ar[d]^{ch^{\rel}_{n,i}} \\
H^{2n-i-1} (X (\C) , \C) / F^n \ar[r]^{\beta} & H^{2n-i-1} (\Gamma (M , \Ah^{\bullet <n})) \; .
}}
\end{equation}
Here, for $n \le i$ the relation of Tamme's relative Chern character
\[
ch^{\rel}_{n,i} : K^{\rel}_i (C^{\infty} (M)) \longrightarrow H^{2n-i-1} (\Gamma (M , \Ah^{\bullet < n}))
\]
to Karoubi's relative Chern character
\[
ch^{\rel}_i : K^{\rel}_i (C^{\infty} (M)) \longrightarrow HC_{i-1} (C^{\infty} (M))
\]
is the following. We have
\begin{align*}
HC_{i-1} (C^{\infty} (M)) & = \Ah^{i-1} (M) / d \Ah^{i-1} (M) \oplus \bigoplus_{k \ge 1} H^{i-1-2k} (M , \C)\\
& = \bigoplus_{k \ge 0} H^{i-1-2k} (\Gamma (M, \Ah^{\bullet < i-k})) \\
& = \bigoplus_{n \le i} H^{2n-i-1} (\Gamma (M , \Ah^{\bullet < n})) \; .
\end{align*}
Tamme's $ch^{\rel}_{n,i}$ is the composition
\begin{equation}
\label{eq:32}
\xymatrix{ch^{\rel}_{n,i} : K^{\rel}_i (C^{\infty} (M)) \ar[r]^-{ch^{\rel}_i} & HC_{i-1} (C^{\infty} (M)) \ar@{->>}[r]^-{\text{projection}} & H^{2n-i-1} (\Gamma (M , \Ah^{\bullet < n}))}
\end{equation}
We are interested in the case $n = i = d+1$. Writing $r_{\Dh}$ for $ch^{\Dh}_{d+1, d+1}$ and noting that $F^{d+1} H^d = 0$, diagrams \eqref{eq:29} and \eqref{eq:31} lead to the following comutative diagrams \eqref{eq:33} and \eqref{eq:34}. The factorization of $\alpha$ via the map $\pi_d$ defined by $\pi_d (z) = \frac{1}{2} (z + (-1)^d \overline{z})$ comes from the exact sequence \eqref{eq:28},
\begin{equation}
\label{eq:33}
\def\objectstyle{\scriptstyle}
\xymatrix{ 
K^{\rel}_{d+1} (X) \ar[d]_{(-1)^d ch^{\rel}_{d+1,d+1}} \ar[rr]^{\theta} & & K_{d+1} (X) \ar[r] \ar[d]^{r_{\Dh}} & K^{\ttop}_{-d-1} (X) \ar[d]^{ch^{\ttop}_{d+1, d+1}} \\
H^d (X (\C) , \C) \ar[r]^-{\pi_d} & H^d (X(\C) , \R (d)) \ar@{^{(}->}[r]^-i & H^{d+1}_{\Dh} (X , \R (d+1)) \ar[r] & H^{d+1} (X (\C) , \R (d+1))
}
\end{equation}
and
\begin{equation}
\label{eq:34}
\vcenter{
\xymatrix{
K^{\rel}_{d+1} (X) \ar[r] \ar[d]_{ch^{\rel}_{d+1,d+1}} & K^{\rel}_{d+1} (C^{\infty} (M)) \ar[d]^{ch^{\rel}_{d+1}} \\
H^d (X (\C), \C) \ar[r]^{\gamma} & HC_d (C^{\infty} (M)) \; .
}}
\end{equation}
Here $\gamma$ sends the class of a smooth closed $d$-form $\omega$ on $X (\C)$ to the class of $\omega \, |_M$ in $\Ah^d (M) / d \Ah^{d-1} (M) \subset HC_d (C^{\infty} (M))$. 

We will also need the following commutative diagram which is induced by the map $R \to C^{\infty} (X (\C)) \to C^{\infty} (M)$.
\begin{equation}
\label{eq:35}
\vcenter{
\xymatrix{
K^{\rel}_{d+1} (X) \ar[r]^{\theta} \ar[d] & K_{d+1} (X) \ar[d] \\
K^{\rel}_{d+1} (C^{\infty} (M)) \ar[r]^{\theta} & K_{d+1} (C^{\infty} (M)) \; .
}}
\end{equation}
We can now finish our conceptual explanation of formula \eqref{eq:2} where we view $m (f)$ as a Fuglede-Kadison determinant. Let $f \in \C [\Z^d] = \C [t^{\pm 1}_1 , \ldots , t^{\pm 1}_d]$ be a Laurent polynomial which does not vanish in any point of $T^d$. By the remark after Corollary \ref{t5.4}, multiplying $f$ by a monomial $t^{\nu_1}_1 \cdots t^{\nu_d}_d$, we may assume that
\[
f \, |_{T^d} \in \Ker (K_1 (C^{\infty} (T^d)) \longrightarrow K^{\ttop}_1 (C^{\infty} (T^d))) \otimes \Q \; .
\]
We have the chain of equalities \eqref{eq:29n} from Corollary \ref{t5.4}
\begin{equation}
\label{eq:36}
m (f) = \log \ddet_{\tau} (f) = \RRe \Delta_{\tau} (f \, |_{T^d}) = \RRe \Delta_{\tau^d} (\{ T^d \} \cdot f \, |_{T^d}) \; .
\end{equation}
Here the cyclic cocyle $\tau^d \in HC^d (C^{\infty} (T^d))$ was the composition 
\begin{equation}
\label{eq:37}
\tau^d : HC_d (C^{\infty} (T^d)) \twoheadrightarrow \Ah^d (T^d) / d \Ah^{d-1} (T^d) \xrightarrow{(2 \pi i)^{-d} \int_{T^d}} \C \; .
\end{equation}
In the algebraic $d$-torus $\Ge^d_{m,\C} = \spec \C [t^{\pm 1}_1 , \ldots , t^{\pm 1}_d]$, consider the complement of the zero-variety of $f$
\[
X = \Ge^d_{m, \C} \setminus V (f) = \spec \C [t^{\pm 1}_1 , \ldots , t^{\pm 1}_d]_f \; .
\]
Since $X$ is an affine variety over $\C$, the complex manifold $X (\C)$ is a Stein manifold and hence $H^i (X (\C) , \Q) = 0$ for $i > d = \dim X$. In particular, the injective map $i$ in diagram \eqref{eq:33} is an isomorphism which we will view as an idenfication. The elements $f , t_1 , \ldots , t_d \in \Oh (X)^{\times} \subset K_1 (X)$ give rise to a symbol
\[
\{ f , t_1 , \ldots , t_d \} = f \cup t_1 \cup \ldots \cup t_d \quad \text{in} \; K_{d+1} (X) \; .
\]
Its images in $K_{d+1} (X) \otimes \Q$ and $K^{\ttop}_{-d-1} (X) \otimes \Q$ lie in the $d+1$-th Adams eigenspaces. Since the topological Chern character $ch^{\ttop}_{d+1 , d+1}$ in diagram \eqref{eq:33} is injective on the $(d+1)$-th Adams eigenspace of $K^{\ttop}_{-d-1} (X) \otimes \Q$ and since the target group $H^{d+1} (X (\C) , \R (d))$ vanishes, it follows that $\{f , t_1 , \ldots , t_d \}$ maps to zero in $K^{\ttop}_{-d-1} (X) \otimes \Q$. Hence by diagram \eqref{eq:33} there is an element $c \in K^{\rel}_{d+1} (X) \otimes \Q$ with $\theta (c) = \{ f , t_1 , \ldots , t_d \}$ and therefore we have the following equality in $H^d (X (\C) , \R (d)) = H^{d+1}_{\Dh} (X , \R (d+1))$
\begin{equation}
\label{eq:38} 
r_{\Dh} (\{ f , t_1 , \ldots , t_d \}) = (-1)^d \pi_d (ch^{\rel}_{d+1 , d+1} (c)) \; .
\end{equation}
Diagram \eqref{eq:34} tells us that the image $c_T \in K^{\rel}_{d+1} (C^{\infty} (T^d)) \otimes \Q$ of $c$ satisfies
\begin{equation}
\label{eq:39}
ch^{\rel}_{d+1} (c_T) = \gamma ch^{\rel}_{d+1,d+1} (c) \quad \text{in} \; HC_d (C^{\infty} (M)) \; .
\end{equation}
Finally from diagram \eqref{eq:35} it follows that we have
\begin{equation}
\label{eq:40}
\theta (c_T) = \{ f \, |_{T^d} , t_1 , \ldots , t_d \} = f \, |_{T^d} \cdot \{ T^d \} = (-1)^d \{ T^d \} \cdot f \, |_{T^d} \; \text{in} \; K_{d+1} (C^{\infty} (T^d)) \; .
\end{equation}
Putting everything together, we find:
\begin{align*}
\langle r_{\Dh} \{ f , t_1 , \ldots , t_d \} , (2 \pi i)^{-d} [T^d] \rangle & \overset{\eqref{eq:38}}{=} (-1)^d \langle \pi_d (ch^{\rel}_{d+1,d+1} (c)) , (2 \pi i)^{-d} [T^d] \rangle \\
& = (-1)^d (2 \pi i)^{-d} \int_{T^d} \pi_d (ch^{\rel}_{d+1,d+1} (c)) \\
& = (-1)^d \RRe (2 \pi i)^{-d} \int_{T^d} ch^{\rel}_{d+1, d+1} (c) \\
& = (-1)^d \RRe \tau^d (\gamma (ch^{\rel}_{d+1,d+1} (c)))\\
& \overset{\eqref{eq:39}}{=} (-1)^d \RRe \tau^d (ch^{\rel}_{d+1} (c_T)) \\
& \overset{\eqref{eq:26}, \eqref{eq:40}}{=} (-1)^d \RRe \Delta_{\tau^d} ((-1)^d \{ T^d \} \cdot f \, |_{T^d}) \\
& = \RRe \Delta_{\tau^d} (\{ T^d \} \cdot f \, |_{T^d} ) \\
& \overset{\eqref{eq:36}}{=} \log \ddet_{\tau} (f) = m (f) \; .
\end{align*}


\end{document}